# A new type of 4D Hybrid Chaos Systems

Reza Parvaz

Department of Mathematics, University of Mohaghegh Ardabili, 56199-11367 Ardabil, Iran.

rparvaz@uma.ac.ir

**Abstract**

In this paper a new type of chaotic system based on sin and logistic systems is introduced. Also the behavior of this new system is studied by using various tests. The results of these tests indicate the appropriate behavior for the proposed new system.

**1-Introduction**

In recent years, the use of chaotic functions has been considered due to its structure. These types of systems are used in image encryption. Also, due to the weaknesses that classical functions have, the use of combination method to improve this group of functions is considered. In this article, we present a four-dimensional chaos system that can be used to create different four-dimensional chaos systems. This system can be used in various algorithms, including the use of image encryption algorithm. It should be noted that all calculations in this article have been done using MATLAB software

**2. 4D Hybrid Chaos Systems**

**2.1 Structure of the Proposed Hybrid System**

In this section, more details about the proposed new hybrid chaos systems based on Tent, Sin and Logistic maps are given. The general structure of the proposed chaos system has been given in Fig. 1. The combination parts in the proposed system based on Tent, Sin and Logistic maps are shown in Fig. 2. The mathematics formulae for each of the parts can be written as follows relations.

First combination part ($A_x$):

$$x_{i+1} = \begin{cases} \alpha_1^x f_1^x \circ F_1^x(r, x_i) + g_1^x(a_1) + h_1^x\left(\dfrac{(\beta_1^x - r)z_i}{2}\right) \bmod 1, & z_i < 0.5, \\ \alpha_2^x f_2^x \circ F_2^x(r, x_i) + g_2^x(a_1) + h_2^x\left(\dfrac{(\beta_2^x - r)(1 - z_i)}{2}\right) \bmod 1, & z_i \geq 0.5. \end{cases}$$

Second combination part ($A_y$):

$$y_{i+1} = \begin{cases} \alpha_1^y f_1^y \circ F_1^y(r, y_i) + g_1^y(a_2) + h_1^y\left(\dfrac{(\beta_1^y - r)\xi_x}{2}\right) \bmod 1, & \xi_x < 0.5, \\ \alpha_2^y f_2^y \circ F_2^y(r, y_i) + g_2^y(a_2) + h_2^y\left(\dfrac{(\beta_2^y - r)(1 - \xi_x)}{2}\right) \bmod 1, & \xi_x \geq 0.5. \end{cases}$$

Third combination part ($A_z$):

$$z_{i+1} = \begin{cases} \alpha_1^z f_1^z \circ F_1^z(r, z_i) + g_1^z(a_3) + h_1^z\left(\dfrac{(\beta_1^z - r)\xi_y}{2}\right) \bmod 1, & \xi_y < 0.5, \\ \alpha_2^z f_2^z \circ F_2^z(r, z_i) + g_2^z(a_3) + h_2^z\left(\dfrac{(\beta_2^z - r)(1 - \xi_y)}{2}\right) \bmod 1, & \xi_y \geq 0.5. \end{cases}$$

Fourth combination part ($A_w$):

$$w_{i+1} = \begin{cases} \alpha_1^w f_1^w \circ F_1^w(r, w_i) + g_1^w(a_4) + h_1^w\left(\dfrac{(\beta_1^w - r)\xi_z}{2}\right) \bmod 1, & \xi_z < 0.5, \\ \alpha_2^w f_2^w \circ F_2^w(r, w_i) + g_2^w(a_4) + h_2^w\left(\dfrac{(\beta_2^w - r)(1 - \xi_z)}{2}\right) \bmod 1, & \xi_z \geq 0.5, \end{cases}$$

Where $a_1 := \{r, x_i, y_i, z_i, w_i\}$, $a_2 = \{r, x_i, x_{i+1}, y_i, z_i, w_i\}$, $a_3 := \{r, x_i, x_{i+1}, y_i, y_{i+1}, z_i, w_i\}$, $a_3 := \{r, x_i, x_{i+1}, y_i, y_{i+1}, z_i, z_{i+1}, w_i\}$, and $\xi_\tau = \tau_i$ or $\tau_{i+1}$, for $\tau = x, y, z$. $F_\varsigma^\tau$, for $\tau = x, y, z, w$, $\varsigma = 1,2$, can be considered as Sin or Logistic maps. $\alpha_\varsigma^\tau, \beta_\varsigma^\tau$, for $\tau = x, y, z, w$, $\varsigma = 1,2$, are arbitrary number in $R$, and $g_\varsigma^\tau, h_\varsigma^\tau$ are considered as arbitrary sufficiently smooth functions. In the proposed system, the best feature of the different chaos maps as Tent, Sin and Logistic maps have been improved by using Composition and transfer operator. In the following, the basic properties of the hybrid system have been studied.

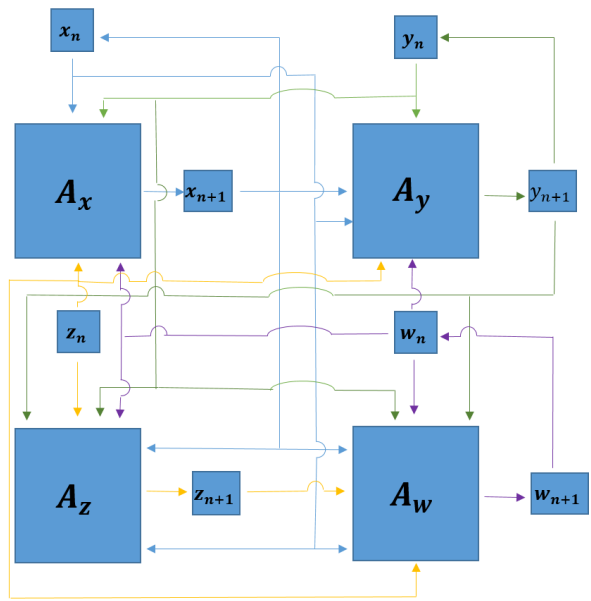

Figure 1. The structure of the proposed chaotic system.

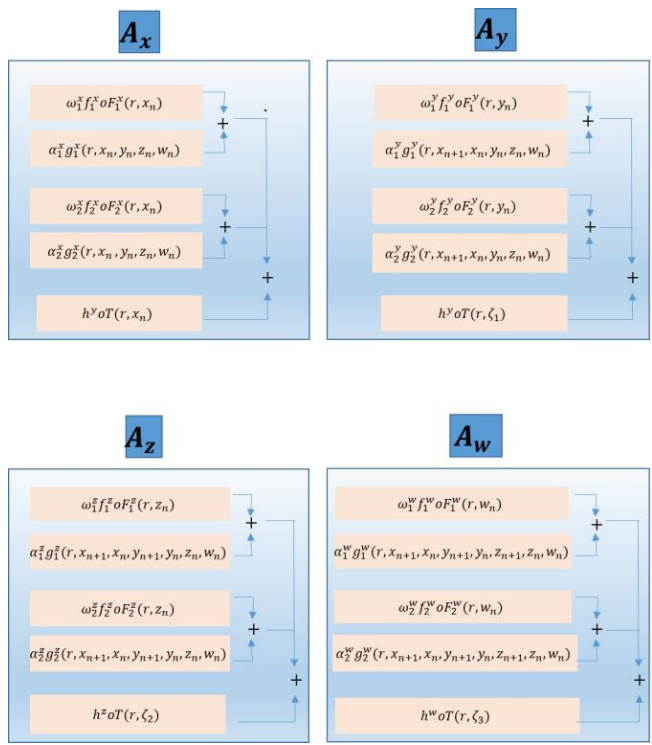

Figure 2. The structure of the combination parts in the proposed system.

In the subsections in order to study hybrid system, the following cases have been considered.

**Case i:** $\{\alpha_1^x,\alpha_2^x,\alpha_1^y,\alpha_2^y,\alpha_1^z,\alpha_2^z,\alpha_1^w,\alpha_2^w\} = \{1,16,10,20,10,20,10,20\}, \{\beta_1^x,\beta_2^x,\beta_1^y,\beta_2^y,\beta_1^z,\beta_2^z,\beta_1^w,\beta_2^w\} = \{6,2,50,30,50,30,50,30\}, \{\xi_x,\xi_y,\xi_z\} = \{x_i,y_i,z_i\}, f_1^x(p) = \cosh(p), f_2^x(p) = \cot(p), f_1^y(p) = f_1^z(p) = f_1^w(p) = p, f_2^y(p) = f_2^w(p) = \sin(\pi p), f_2^z(p) = \exp(\pi p),\ g_1^x(a_1) = 15\tanh(rx_i + z_i) + \sin(w_i) + 12\cos(rx_i),\ g_2^x(a_1) = -7ry_i + \exp(1 + 2w_i) + z_i + 7\log(\pi rx_i),\ g_1^y(a_2) = 2\tan(rx_i + y_i + 2z_i + w_i),\ g_2^y(a_2) = z_i + w_i + 14\exp(20rx_i),\ g_1^z(a_3) = 2\tan(rx_i + y_i) + w_i + z_i,\ g_2^z(a_3) = 14\exp(20rx_i + w_i) + \sin(z_i),\ g_1^w(a_4) = 2\tan(rx_i + y_i + z_i) + w_i,\ g_2^w(a_4) = 14\exp(20rx_i + w_i) + z_i, h_1^x(p) = \sin(2p),\ h_2^x(p) = 4p, h_2^y(p) = \cot(p),\ h_1^y(p) = h_1^z(p) = h_1^w(p) = \exp(2p), h_2^z(p) = h_2^w(p) = \cot(4p)$.

**Case ii:** $\{\alpha_1^x,\alpha_2^x,\alpha_1^y,\alpha_2^y,\alpha_1^z,\alpha_2^z,\alpha_1^w,\alpha_2^w\} = \{7,2,4,4,3,5,5,5\}, \{\beta_1^x,\beta_2^x,\beta_1^y,\beta_2^y,\beta_1^z,\beta_2^z,\beta_1^w,\beta_2^w\} = \{1,2,3,3,1,2,2,2\}, \{\xi_x,\xi_y,\xi_z\} = \{x_{i+1},y_{i+1},z_{i+1}\}, f_1^x(p) = f_2^w(p) = \cos(p), f_2^x(p) = f_2^y(p) = f_2^z(p) = f_1^w(p) = p, f_1^z(p) = \exp(p),\ g_1^x(a_1) = \sin(w_i + z_i) + rx_iy_i, g_2^x(a_1) = \sin(ry_i + x_i) + \log(7 + w_i + z_i),\ g_1^y(a_2) = rx_i + y_i + \exp(rx_{i+1}) + \cos(z_i + w_i), g_2^y(a_2) = z_i - w_i + \log(20rx_{i+1} + x_i), g_1^z(a_3) = \cot(rx_{i+1} + y_{i+1}) + \sin(x_i + w_iz_i), g_2^z(a_3) = \exp(x_i + w_i + z_i) + \sin(x_{i+1} + y_{i+1}), g_1^w(a_4) = 2\cot(rx_{i+1} + y_{i+1} + z_{i+1}) + \log(x_i + w_i), g_2^w(a_4) = \exp(ry_{i+1} + x_{i+1} + 2w_i) + y_i + z_i, h_1^x(p) = \cos(20p),\ h_2^x(p) = 5p, h_1^y(p) = \log(4p), h_2^y(p) = \cos(6p),\ h_1^z(p) = h_1^w(p) = \exp(4p), h_2^z(p) = h_2^w(p) = \cos(p)$.

Also, in the case i, $\{F_1^x, F_1^z, F_2^z\}$ are considered as Tent map, and $\{F_2^x, F_1^y, F_2^y, F_2^z\}$ are considered as Sin map. For case ii, $\{F_1^x, F_1^z\}$ and $\{F_2^x, F_1^y, F_2^y\}$ are considered as Tent and Sin maps, respectively.

## 2.1 Chaotic Behavior Analysis

In this subsection same important tests for the proposed chaos system are discussed. One of the important values in the study of the behavior of the chaos system is Lyapunov exponent or Lyapunov characteristic exponent. A $n$-dimensional chaos systems in general have $n$ values for Lyapunov exponent. There are many methods for calculating this value as [1-3]. The method based on QR algorithm has been used for obtained Lyaponov exponent in Fig. 3 for the case i. More details about this method can be found in [3]. The positive or negative values of the resulting values are related with the structure of a chaos system. This relation had been studied in many papers. In [4], the relation has been given as follows '' *In an n-dimensional dynamical system we have n Lyapunov exponents. Each $\lambda_k$ represents the divergence of k-volume. The sign of the Lyapunov exponents indicates the behavior of nearby trajectories. A negative exponent indicates that neighboring trajectories converge to the same trajectory A positive exponent indicates that neighboring trajectories diverge* [4]''. Also the following theorem has been given in [5] for this value.

**Theorem 2.1** *If at least one of the average Lyapunov exponents is positive, then the system is chaotic; if the average Lyapunov exponent is negative, then the orbit is periodic and when the average Lyapunov exponent is zero, a bifurcation occurs.*

The results in Fig. 3 show that all four value of Lyapunov exponent in the proposed system are positive. Then by using above studies, we can say that in the proposed system in the all neighboring trajectories diverge. In order to compare proposed system, the Lyapunov exponent has been compared with 4D Chaotic Laser System [6] in Fig. 3. It is observed that the proposed system has better chaos behavior then Chaotic Laser System.

Another tool for study chaotic behavior is bifurcation analysis. In the Fig. 4, the results for the bifurcation analysis of the case ii of the proposed system have been shown. The chaotic attractors can be studied by this figure. The attractor for $r \in (0, 1.2]$ is given in the vertical line at that $r$. Also, the cobweb plot (or Verhulst diagram) for case i have been given in the Fig.s 5. By using this results, it is observed that for the given values, the resulting sequences of the proposed system has chaotic behavior.

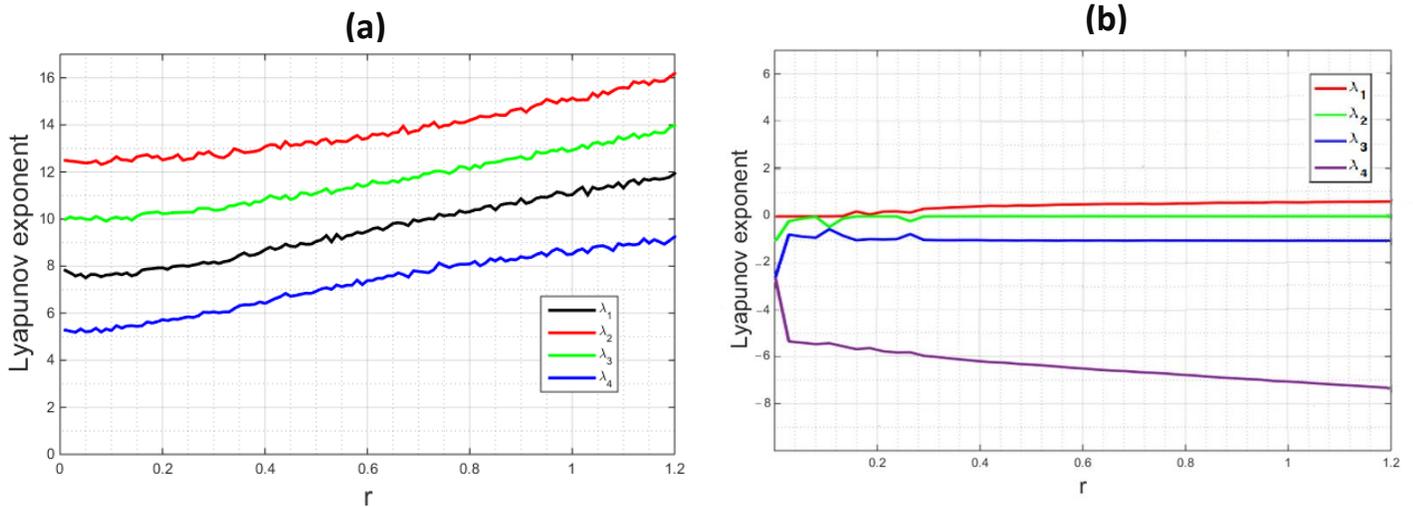

Figure 3. Lyapunov exponent values for: (a) Case i, (b) Chaotic system in [6].

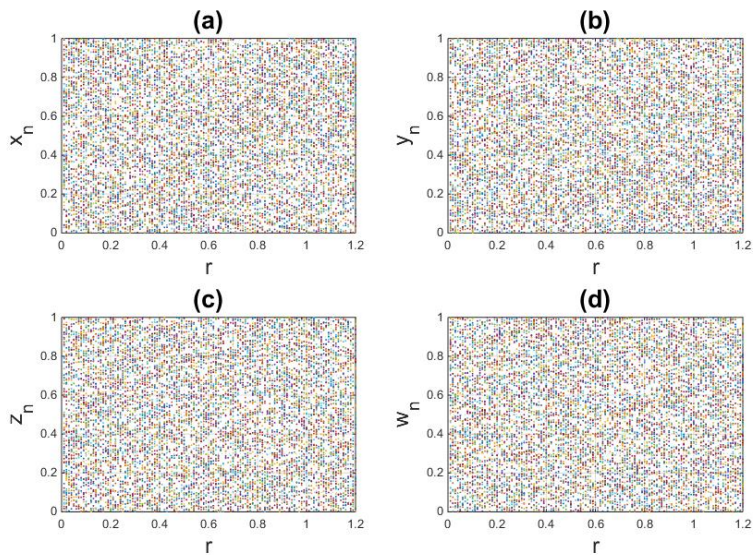

Figure 4. Bifurcation diagram results for the case ii.

Distribution is another important factor in chaotic system. One of the reasons for the weakness of the statistical attack is nonuniform distribution. The histogram plots of the proposed system for the case i are given in the Fig. 6. Also, the distribution patterns of the case ii are shown in Fig. 7. By using these results, it can be seen that the generated sequence of the proposed maps have a flat distribution.

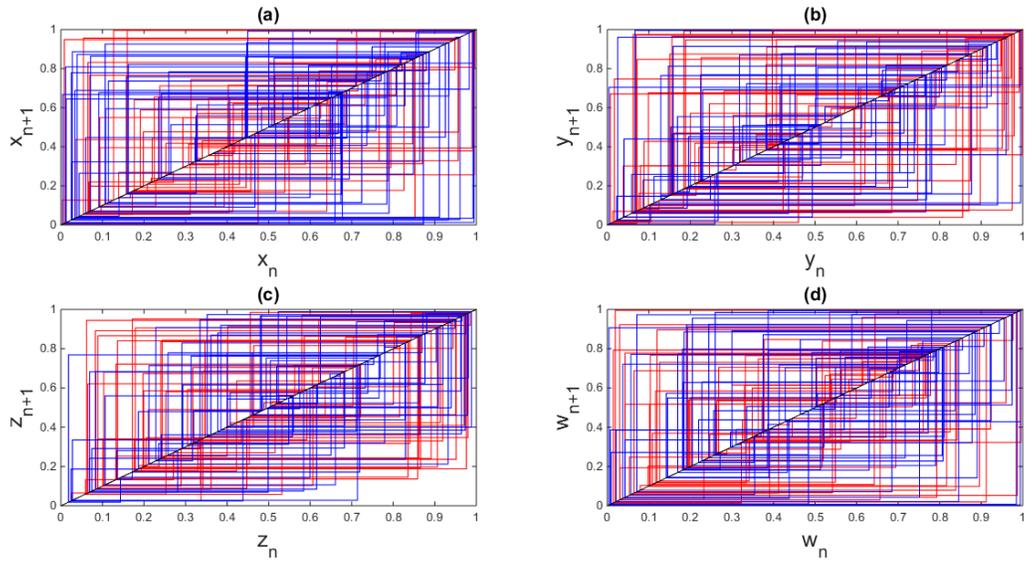

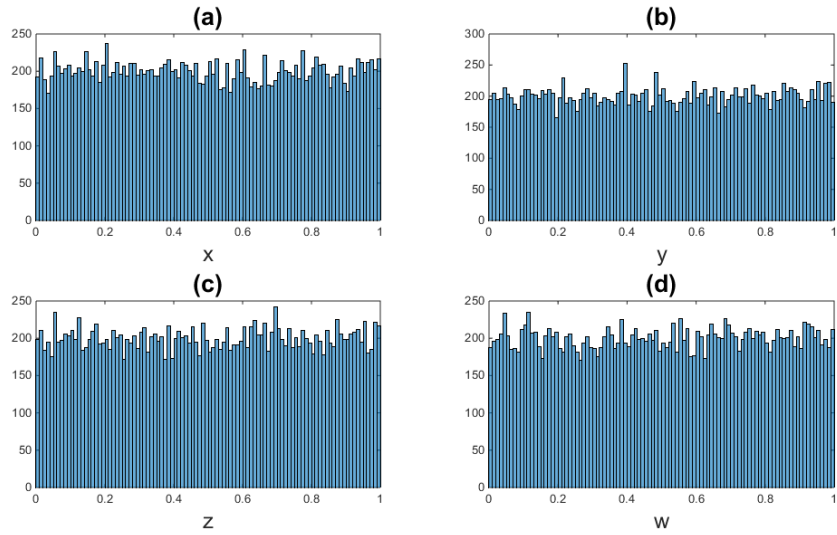

Figure 6. Histogram plots of the case i for $r = 0.5$.

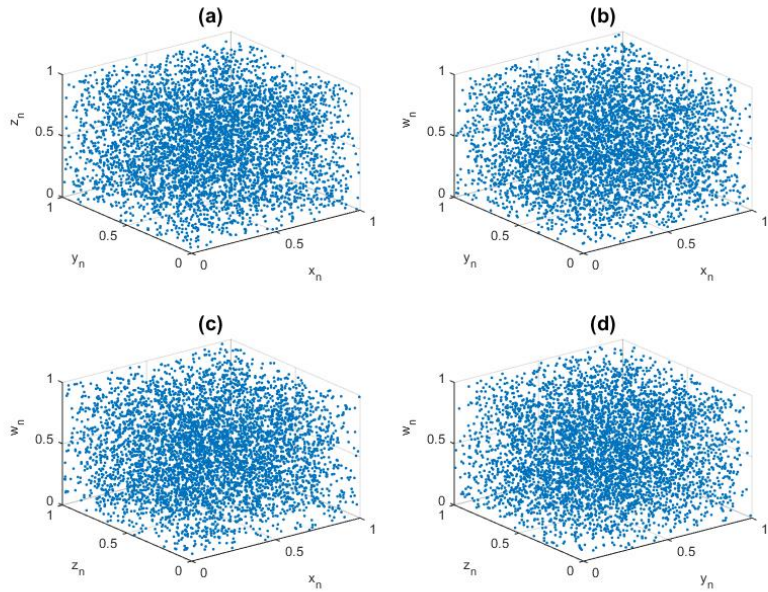

Figure 6. Distribution patterns of the case ii for $r = 0.4$.